\documentclass[11pt,a4paper]{article}
\usepackage{amsmath}
\usepackage{amsfonts}
\usepackage{amssymb}
\usepackage{amsthm}
\usepackage{amscd}
\usepackage{graphicx}
\theoremstyle{plain}

\theoremstyle{definition}

\newtheorem{remark}{Remark}[section]

\title{Recent developments of biharmonic conjecture and modified biharmonic conjectures}

\author{
	\emph{Bang-Yen Chen}
	 \\
	Michigan State University \\
	Department of Mathematics \\ 
	619 Red Cedar Road \\ East Lansing, Michigan 48824, USA
	\\e-mail: \emph{bychen@math.msu.edu}
     }

\date{}

\begin{document}

\maketitle

\begin{abstract}
A submanifold $M$ of a Euclidean $m$-space $\mathbb E^m$ is said to be biharmonic if
$\Delta H=0$ holds identically, where $ H$ is the mean curvature vector field and $\Delta$ is the Laplacian on $M$. In 1991, the author conjectured in  \cite{C91} that every biharmonic submanifold of a Euclidean space is minimal. The study of  biharmonic submanifolds is nowadays a very active research subject. In particular, since 2000 biharmonic submanifolds have been
receiving a growing attention and have become a popular subject of study with
many progresses. 

In this article, we provide a brief survey on recent developments concerning my original conjecture and  generalized biharmonic conjectures. At the end of this article, I present two modified conjectures related with original biharmonic conjecture.

\end{abstract}

\section{Introduction}
Let $x:M\to {\mathbb E}^m$ be an isometric immersion from a Riemannian $n$-manifold into a Euclidean $m$-space. 
Denote by $\Delta, x$ and $H$ the Laplacian, the position vector  and the mean curvature vector  of $M$, respectively. 
Then $M$ is called a biharmonic submanifold if $\Delta  H=0$. Due to the well-known Beltrami's formula, $\Delta  x=-n H$, 
it is obvious that every minimal submanifold of $\mathbb E^m$ is a biharmonic submanifold.

The study of biharmonic submanifolds was initiated by the author in the middle of 1980s in his program of understanding the finite type submanifolds in Euclidean spaces; also independently by G. Y. Jiang  \cite{J} for his study of Euler-Lagrange's equation of bienergy functional in the sense of Eells and Lemaire. 

The author showed in 1985 that biharmonic surfaces in $\mathbb E^3$ are minimal (unpublished then, also independently by Jiang \cite{J}). This result was the starting point of I. Dimitric's work on his doctoral thesis at Michigan State University (cf. \cite{Dim89}).  
In particular, Dimitric extended author's unpublished result to  biharmonic hypersurfaces of $\mathbb E^m$ with at most  two distinct principal curvatures \cite{Dim89}.  In his thesis,
Dimitric also proved that every biharmonic submanifold of finite type 
in ${\mathbb{E}}^m$ is minimal.
Another extension of this result on biharmonic surfaces was given by T. Hasanis and T. Vlachos in \cite{HV} (see also
 \cite{D98}). They proved that biharmonic hypersurfaces of ${\mathbb{E}}^4$ are minimal.

Formally, the author made in \cite{C91} the following.
\vskip.06in

\hskip.3in {\bf Biharmonic Conjecture}:
\emph{The only biharmonic submanifolds of Euclidean spaces are 
the minimal ones.}
\vskip.06in

A {\it biharmonic map} is a map $\phi:(M,g)\to (N,h)$ between Riemannian manifolds that is a critical point of the bienergy functional:
\begin{align}E^2(\phi,D)=\frac{1}{2}\int_D ||\tau_\phi||^2* 1\end{align}
for every compact subset $D$ of $M$, where
$\tau_\phi={\rm trace}_g\nabla d\phi$ is the tension field $\phi$. 
  The Euler-Lagrange equation of this functional gives the biharmonic map
equation (see \cite{J})
\begin{align}\label{BH} \tau^2_\phi:={\rm trace}_g(\nabla^\phi\nabla^\phi-\nabla^\phi_{\nabla^M})\tau_\phi-{\rm trace}_g R^N(d\phi,\tau_\phi)d\phi=0,\end{align}
where $R^N$ is the curvature tensor of $(N,h)$.
Equation \eqref{BH} states that $\phi$ is a biharmonic map if and only if its bi-tension field $\tau^2_\phi$ vanishes.  

Let $M$ be an $n$-dimensional submanifold of a Euclidean $m$-space $\mathbb E^m$.
 If we denote by $\iota:M\to \mathbb E^m$ the inclusion map of the submanifold, then the tension field of the inclusion map is given by $\tau_\iota=-\Delta\iota=-n H$ according to Beltrami's formula. Thus $M$ is a biharmonic submanifold if and only if 
$$n\Delta  H=- \Delta^2 \iota=-\tau^2_\iota=0,$$
i.e., the inclusion map $\iota$ is a biharmonic map.

 Caddeo,  Montaldo and  Oniciuc \cite{CMO02} proved that every biharmonic surface in 
the hyperbolic $3$-space $H^3(-1)$ of constant curvature $-1$ is minimal. They also proved that biharmonic hypersurfaces of
$H^n(-1)$ with at most two distinct principal curvatures are minimal \cite{CMO01}. 

Based on these, Caddeo, Montaldo and Oniciuc made in  \cite{CMO01} the following.
\vskip.06in

\hskip.3in  {\bf The generalized Chen's conjecture:}
\emph{Any biharmonic submanifold of a Riemannian manifold with non-positive sectional curvature 
 is minimal.}
\vskip.06in

The study of  biharmonic submanifolds is nowadays a very active research subject. In particular, since 2000 biharmonic submanifolds have been
receiving a growing attention and have become a popular subject of study with
many progresses. 

 In this article, we provide a brief survey on recent developments concerning my original conjecture and  generalized biharmonic conjectures. At the end of this article, I present two modified conjectures related with original biharmonic conjecture.

\section{Recent developments  on Chen's original biharmonic conjecture}

 Let  $x :M\to \mathbb E^m$ be an isometric immersion of a Riemannian $n$-manifold $M$ into a Euclidean $m$-space $\mathbb E^m$.  
 Then $M$ is biharmonic if and only if it satisfies the following fourth order strongly elliptic semi-linear PDE system (see, for instance, \cite{C84,C96,C11})
 \begin{align}\notag &\begin{cases}\; \Delta^{D} H + \sum_{i=1}^{n}
\sigma(A_{H}e_{i},e_{i})=0,\\  \\ \;n\,\nabla\!\left<\right.\! H, H\! \left.\right> + 4\, {\rm trace}\,A_{D H}=0,\end{cases}\end{align}
where $\Delta^D$ is the Laplace operator associated with the normal connection $D$, $\sigma$ the second fundamental form, $A$ the shape operator, $\nabla\!\left<\right.\! H, H\! \left.\right>$ the gradient of the squared mean curvature, and $\{e_1,\ldots,e_n\}$ an orthonormal frame of $M$.

An immersed submanifold $M$ in a Riemannian manifold $N$ is said to be {\it properly
immersed} if the immersion is a proper map, i.e.,  the preimage of each compact set in
$N$ is compact in $M$.

The {\it total mean curvature} of a submanifold $M$ in a Riemannian manifold is given by $\int_M | H|^2 dv$.

 Denote by $K(\pi)$ the sectional curvature of a given Riemannian $n$-manifold $M$ associated with a plane section $\pi\subset T_pM$, $p\in M$. For any orthonormal basis $e_1,\ldots,e_n$ of the tangent space $T_pM$, the scalar curvature $\tau$ at $p$ is defined to be $\tau(p)=\sum_{i<j} K(e_i\wedge e_j). $

Let $L$ be a subspace of $T_pM$  of dimension $r\geq 2$  and $\{e_1,\ldots,e_r\}$ an orthonormal basis of $L$. The scalar curvature $\tau(L)$ of $L$ is defined by
$$\tau(L)=\sum_{\alpha<\beta} K(e_\alpha\wedge e_\beta),\quad 1\leq \alpha,\beta\leq r.$$

For an integer  $r\in [2,n-1]$, the {\it $\delta$-invariant} $\delta(r)$ of $M$ is defined  by  (cf. \cite{C00,C11})
\begin{align}\label{1.3} \delta(r)(p)=\tau(p)- \inf\{\tau(L)\},\end{align} where $L$ run over all $r$-dimensional linear subspaces of $T_pM$. 

For any $n$-dimensional submanifold $M$ in $\mathbb E^m$ and any integer $r\in [2, n-1]$,  the author proved the following general sharp inequality  (cf. \cite{C00,C11}):
\begin{align}\label{1.4} \delta(r) \leq  \frac{n^2(n-r)}{2(n-r+1)} |\overrightarrow H|^2.\end{align}

 A submanifold  in $\mathbb E^m$ is called {\it $\delta(r)$-ideal} if it satisfies the equality case of \eqref{1.4} identically. Roughly speaking ideal submanifolds are submanifolds which receive the least possible tension from its ambient space (cf. \cite{C00,C11}).
 
 A hypersurface of a Euclidean space is called {\it weakly convex} if it has non-negative principle curvatures.
 
It follows immediately from the definition of biharmonic submanifolds and Hopf's lemma that every biharmonic submanifold in a Euclidean space is non-compact.

\vskip.06in
The following provides an overview of some affirmative partial solutions to my original biharmonic conjecture.

\begin{itemize}
\item Biharmonic surfaces in $\mathbb E^3$ (B.-Y. Chen \cite{C91,C11} and G. Y. Jiang \cite{J}).

\item Biharmonic  curves (I. Dimitric \cite{Dim89,Dim92}).

\item Biharmonic  hypersurfaces in $\mathbb E^4$ (T. Hasanis and T. Vlachosin \cite{HV}) (a different proof by F. Defever  \cite{D98}).

\item Spherical submanifolds (B.-Y. Chen \cite{C96}).

\item Biharmonic hypersurfaces with at most 2 distinct principle curvatures (I. Dimitric  \cite{Dim89}).

\item Biharmonic  submanifolds of finite type (I. Dimitric \cite{Dim89,Dim92}).

\item Pseudo-umbilical biharmonic  submanifolds (I. Dimitric \cite{Dim92}).

\item Biharmonic  submanifolds which are complete and proper (Akutagawa and Maeta  \cite{AM}).

\item Biharmonic properly immersed submanifolds (S. Maeta \cite{M12a}).

\item Biharmonic  submanifolds  satisfying the
decay condition at infinity $$\lim_{\rho\to \infty}\frac{1}{\rho^2}\int_{f^{-1}(B_\rho)}| H |^{2}dv=0,$$
where $f$ is the immersion, $B_\rho$ is a geodesic ball of $N$ with radius $\rho$  (G. Wheeler \cite{Wh}).

\item Submanifolds whose $L^p,\, p\geq 2$,  integral of the mean curvature vector field satisfies
certain decay condition at infinity (Y. Luo  \cite{Luo3}).

\item $\delta(2)$-ideal and $\delta(3)$-ideal biharmonic  hypersurfaces (B.-Y. Chen and M. I. Munteanu  \cite{CM}).

\item Weakly convex biharmonic submanifolds (Y. Luo in \cite{Luo1}).
\end{itemize}

In \cite{Ou09}, Y.-L. Ou constructed examples to show that my original biharmonic conjecture cannot be generalized to the case of biharmonic conformal submanifolds in Euclidean spaces. 
\vskip.05in

\begin{remark} {\sc My original biharmonic conjecture is  still open}.
\end{remark}

\begin{remark} My biharmonic conjecture is false if the ambient Euclidean space were replaced by a pseudo-Euclidean space. The simplest examples are constructed by Chen and Ishikawa in \cite{CI91}. For instance, we have the following.
\vskip.05in

{\bf Example.}  Let $f(u,v)$ be a proper biharmonic function, i.e. $\Delta f\ne 0$ and $\Delta^2 f=0$. Then 
\begin{align}\label{Bi} x(u,v)=(f(u,v),f(u,v),u,v)\end{align}
defines a biharmonic, marginally trapped surface in the Minkowski 4-space $\mathbb E^4_1$ endowed with the Lorentzian metric $g_0=-dt_1^2+dx_1^2+dx_2^2+dx_3^2$. 

Here, by a marginally trapped surface, we mean a space-like surface in $\mathbb E^4_1$ with light-like mean curvature vector field.

It was proved in \cite{CI91} that the biharmonic surfaces defined by \eqref{Bi} are the only biharmonic, marginally trapped surfaces in $\mathbb E^4_1$.
\end{remark}

\section{Recent developments  on Caddeo-Montaldo-Oniciuc's Generalized Chen's biharmonic conjecture}

Let $M$ be a submanifold of a Riemannian manifold with inner product $\left<\;\,,\;\right>$, then $M$ is called {\it $\epsilon$-superbiharmonic} if
$$\left<\right.\!\Delta  H, H\! \left.\right>\geq (\epsilon-1)|\nabla H|^2,$$
where $\epsilon\in [0,1]$ is a constant.
For a complete Riemannian manifold $(N,h)$ and $\alpha\geq 0$, if the sectional curvature $K^N$ of $N$ satisfies
$$K^N\geq -L(1+{\rm dist}_N(\,\cdot\,,q_0)^2)^{\frac{\alpha}{2}}$$
for some $L>0$ and $q_0\in N$, 
then we  call that $K^N$ has a polynomial growth bound of order $\alpha$ from below.
\vskip.06in

There are also many affirmative partial answers to the generalized Chen's biharmonic conjecture. The following provides a brief overview of the affirmative partial answers to this generalized conjecture.

\begin{itemize}
\item Biharmonic hypersurfaces in the hyperbolic 3-space $H^3(-1)$ (Caddeo,  Montaldo and  Oniciuc \cite{CMO01}).

\item  Biharmonic hypersurfaces in $H^4(-1)$ (Balmu\c s, Montaldo and  Oniciuc \cite{BMO10b}).

\item Pseudo-umbilical biharmonic submanifolds of $H^m(-1)$ (Caddeo,  Montaldo and  Oniciuc \cite{CMO01}).

\item Biharmonic hypersurfaces of $H^{n+1}(-1)$ with at most two distinct principal curvatures (Balmu\c s, Montaldo and  Oniciuc \cite{BMO08}).

\item Totally umbilical biharmonic hypersurfaces in Einstein spaces (Y.-L. Ou \cite{Ou10}).

\item Biharmonic hypersurfaces with finite total mean curvature in a Riemannian manifold of non-positive Ricci curvature (Nakauchi and Urakawa \cite{NU1}).

\item Biharmonic submanifolds with finite total mean curvature in a Riemannian manifold of non-positive sectional curvature (Nakauchi and Urakawa \cite{NU2}).

\item Complete biharmonic hypersurfaces $M$ in a Riemannian manifold of non-positive Ricci curvature whose mean curvature vector satisfies $\int_M | H|^\alpha dv<\infty$ for some $\epsilon>0$ with $1+\epsilon\leq \alpha<\infty$ (S. Maeta \cite{M13}).

\item Biharmonic properly immersed submanifolds in a complete Riemannian manifold with non-positive sectional curvature whose sectional curvature has polynomial growth bound of order less than 2 from below (S. Maeta \cite{M12b}).

\item Complete biharmonic submanifolds  with finite bi-energy and energy in a non-positively curved Riemannian manifold (N. Nakauchi, H. Urakawa and S. Gudmundsson \cite{NUG}).

\item Complete oriented biharmonic hypersurfaces $M$  whose mean curvature $H$ satisfying $H\in L^2(M)$  in a Riemannian manifold with non-positive Ricci tensor (Al\'{\i}as, Garc\'{\i}a-Mart\'{\i}nez and Rigoli \cite{Al}).

\item Compact biharmonic submanifolds in a Riemannian manifold with non-positive sectional curvature (G.-Y. Jiang \cite{J} and S. Maeta \cite{M13}).

\item $\epsilon$-superbiharmonic submanifolds in a complete Riemannian manifolds satisfying the decay condition at infinity
$$\lim_{\rho\to \infty}\frac{1}{\rho^2}\int_{f^{-1}(B_\rho)}| H |^{2}dv=0,$$
where $f$ is the immersion, $B_\rho$ is a geodesic ball of $N$ with radius $\rho$  (G. Wheeler \cite{Wh}).

\item Complete biharmonic submanifolds (resp., hypersurfaces) $M$ in a Riemannian manifold of non-positive sectional (resp., Ricci) curvature whose mean curvature vector satisfies $\int_M | H^p |dv<\infty$ for some $p>0$  (Y. Luo \cite{Luo2}).

\item Complete biharmonic submanifolds (resp., hypersurfaces) in a Riemannian manifold
whose sectional curvature (resp., Ricci curvature) is non-positive with at most polynomial
volume growth (Y. Luo \cite{Luo2}).

\item Complete biharmonic submanifolds (resp., hypersurfaces) in a negatively curved Riemannian
manifold whose sectional curvature (resp., Ricci curvature) is smaller that $-\epsilon$ for some
$\epsilon>0$ (Y. Luo \cite{Luo2}).

\item Proper $\epsilon$-superharmonic submanifolds $M$ with $\epsilon>0$ in a complete Riemannian manifold $N$ whose mean curvature vector satisfying the growth condition
$$\lim_{\rho\to \infty}\frac{1}{\rho^2}\int_{f^{-1}(B_\rho)}| H |^{2+a}dv=0,$$
where $f$ is the immersion, $B_\rho$ is a geodesic ball of $N$ with radius $\rho$, and $a\geq 0$ (Luo \cite{Luo2}).

\end{itemize}

On the other hand, it was proved by Y.-L. Ou and L. Tang in \cite{OT}  that the generalized Chen's biharmonic conjecture is false in general by constructing foliations of proper biharmonic hyperplanes in a $5$-dimensional conformally flat space with negative sectional curvature.

Further counter-examples were constructed  in \cite{LO} by T. Liang and Y.-L. Ou.

\section{Maeta's Generalized Chen's conjecture}

A submanifold of a Euclidean space is called {\it $k$-harmonic} if its mean curvature vector satisfies $\Delta^{k-1} \overrightarrow H=0$.  
It follows from Hopf's lemma that such submanifolds are always non-compact.
Some relationships between $k$-harmonic and harmonic maps  of Riemannian manifolds into Euclidean $m$-space $\mathbb E^{m}$ have been obtained by the author in \cite{C07}. 

Recently, S. Maeta \cite{M12} found  some relations between $k$-harmonic and harmonic maps of Riemannian manifolds into 
non-flat real space forms.

It follows from   \cite[Proposition 3.1]{C07} that every $k$-harmonic submanifold of $\mathbb E^{m}$ is either minimal or of infinite type (in the sense of \cite{C84}). 
On the other hand, it is also well-known that all $k$-harmonic curves in $\mathbb E^{m}$ are of finite type  (see \cite[Proposition 4.1]{CP}). 
Consequently, every $k$-harmonic curve in $\mathbb E^{m}$ is an open portion of line (this known fact was rediscovered recently in \cite[Theorem 5.5]{M12}). 

Based on this fact, S. Maeta  \cite{M12} made another generalized Chen's conjecture; namely,

\vskip.05in

\hskip.3in
``{\it The only $k$-harmonic submanifolds of a Euclidean space are the minimal  ones.}''

\section{Two related biharmonic conjectures}

Finally, I present two  biharmonic conjectures related to my original biharmonic conjecture.
\vskip.06in

\hskip.3in {\bf Biharmonic Conjecture for Hypersurfaces}:
\emph{Every biharmonic hypersurface of Euclidean spaces is minimal.}
\vskip.06in

The global version of my original biharmonic conjecture can be found, for instance, in \cite{AM,M13}.
\vskip.06in

\hskip.3in {\bf Global Version of Chen's biharmonic Conjecture}:
\emph{Every complete biharmonic submanifold of a Euclidean space is minimal.}

\end{document}